\documentclass[a4paper,11pt,reqno, english]{amsart}

\usepackage{amsmath,amssymb,amscd,amsthm,amsfonts}
\usepackage{graphicx,subfigure}
\usepackage{hyperref}
\usepackage{dsfont}
\usepackage{color}
\usepackage{paralist}

\usepackage[left=1.2in,right=1.2in,top=1in,bottom=1.05in]{geometry}

\newtheorem{theorem}{Theorem}[section]

\newtheorem{lemma}[theorem]{Lemma}
\newtheorem{claim}[theorem]{Claim}
\newtheorem{corollary}[theorem]{Corollary}
\newtheorem{conjecture}[theorem]{Conjecture}
\newtheorem{remark}[theorem]{Remark}

\newcommand{\rr}{\mathds{R}}

\newcommand{\Z}{\mathbb{Z}}

\newcommand{\R}{\mathds{R}}

\DeclareMathOperator{\conn}{conn}

\title[The topological Tverberg problem]{The topological Tverberg problem \\ beyond prime powers}

\author[Frick]{Florian Frick}
\address{Dept.\ Math.\ Sciences, Carnegie Mellon University, Pittsburgh, PA 15213, USA}
\email{frick@cmu.edu} 

\author[Sober\'on]{Pablo Sober\'on}\address{Baruch College, City University of New York, One Bernard Baruch Way, New York, NY 10010, USA} 
\email{pablo.soberon-bravo@baruch.cuny.edu}

\thanks{FF was supported by NSF grant DMS 1855591 and a Sloan Research Fellowship. PS was supported by NSF grant DMS 1851420 and PSC-CUNY grant 62639-00-50.}

\date{May 11, 2020, revised August 1, 2023}

\begin{document}

\maketitle

\begin{abstract}
\small
	Tverberg-type theory aims to establish sufficient conditions for a simplicial complex $\Sigma$ such that every continuous map $f\colon \Sigma \to \R^d$ maps $q$ points from pairwise disjoint faces to the same point in~$\R^d$. Such results are plentiful for $q$ a power of a prime. However, for $q$ with at least two distinct prime divisors, results that guarantee the existence of $q$-fold points of coincidence are non-existent---aside from immediate corollaries of the prime power case. Here we present a general method that yields such results beyond the case of prime powers. In particular, we prove previously conjectured upper bounds for the topological Tverberg problem for all~$q$.
\end{abstract}

\section*{Erratum}

We are grateful to an anonymous referee, who has pointed out a mistake in our proof. The error is in Claim 5.4, where the induction fails in the case $k=q$. The proof of our main result Theorem~1.1 relies on Claim 5.4. Thus the proof of Theorem~1.1 is incomplete. The reduction in Section~4 is unaffected. We apologize for this mistake.

\section{Introduction}

In 1959 Bryan Birch~\cite{birch1959} proved that \textit{in any straight-line drawing of the complete graph $K_{3q}$ on $3q$ vertices there are $q$ pairwise vertex-disjoint $3$-cycles that surround a common point}. Equivalently, for $3q$ points in the plane, there is a partition into $q$ sets whose convex hulls all share a common point. In this phrasing of the result, one may observe that $3q-2$ points suffice for such a partition to exist. Helge Tverberg~\cite{tverberg1966} generalized this result to higher dimensions: \textit{Any $(q-1)(d+1)+1$ points in $\R^d$ may be partitioned into $q$ sets whose convex hulls all share a common point}.

One may also wonder whether it is necessary that the drawing of $K_{3q}$ is a straight-line drawing in Birch's result. Indeed, Imre B\'ar\'any conjectured a topological generalization of Tverberg's theorem: \textit{Any continuous map $f\colon \Delta_{(q-1)(d+1)} \to \R^d$ from the $(q-1)(d+1)$-dimensional simplex~$\Delta_{(q-1)(d+1)}$ to $\R^d$ identifies points from $q$ pairwise disjoint faces}. For an affine map $f$ this is Tverberg's theorem since the affine image of a face is the convex hull of (the images of) its vertices. B\'ar\'any's topological Tverberg conjecture was proven for $q$ a prime by B\'ar\'any, Shlosman, and Sz\H ucs~\cite{barany1981} and for $q$ a prime power by \"Ozaydin~\cite{ozaydin1987} and Volovikov~\cite{volovikov1996}. Recently, counterexamples were exhibited for any $q$ that is not a power of a prime; see~\cite{blagojevic2019, frick2015} and Mabillard and Wagner~\cite{mabillard2014, mabillard2015}. Tverberg-type problems have been of significant interest; see~\cite{barany2016, barany2018, blagojevic2017, deloera2019} for recent surveys.

These developments leave open the problem of topological generalizations of Birch's original result and its higher-dimensional versions. There are no non-trivial upper bounds for the topological Tverberg problem beyond the case of prime powers. The purpose of the present manuscript is to prove a continuous generalization of Birch's result in any dimension and for all integers~${q\ge2}$.

\begin{theorem}
\label{thm:main}
	Let $q\ge2$ and $d\ge1$ be integers. Let $f\colon \Delta_{q(d+1)-1} \to \R^d$ be a continuous map. Then there are $q$ pairwise disjoint faces $\sigma_1, \dots, \sigma_q$ of $\Delta_{q(d+1)-1}$ such that $f(\sigma_1) \cap \dots \cap f(\sigma_q) \ne \emptyset$.
\end{theorem}

Avvakumov, Karasev, and Skopenkov~\cite{avvakumov2019} show that if $q$ is not a prime power then there exists a continuous map $f\colon \Delta_n \to \R^d$ with $f(\sigma_1) \cap \dots \cap f(\sigma_q) = \emptyset$ for any $q$ pairwise disjoint faces $\sigma_1, \dots, \sigma_q$ for $n = q\left(d+1 - \Big\lceil \frac{d+2}{q+1}\Big \rceil\right) -2$, which is currently the best lower bound for the topological Tverberg problem for large~$d$. Theorem~\ref{thm:main} shows that the multiplicative factor for the dimension of the simplex is asymptotically at most~$q$ for $q$ not a power of a prime.

For $q$ a prime power, Theorem~\ref{thm:main} is weaker than the topological generalization of Tverberg's theorem.  If $q+1$ is a prime power, it is a simple consequence thereof; see Section~\ref{sec:simple}. Blagojevi\'c, the first author, and Ziegler~{\cite[Conj.~5.5]{blagojevic2019}} conjectured Theorem~\ref{thm:main}.  They also conjectured that the bound is optimal, which remains open. The best bounds on $n$ such that any continuous map $f \colon \Delta_n \to \R^q$ exhibits a $q$-fold point of coincidence among pairwise vertex-disjoint faces are simply derived by choosing $n$ sufficiently large to guarantee such an intersection for $p$ faces, where $p\ge q$ is a prime power. Indeed, the topological tools used to prove these results---the non-existence of associated equivariant maps---fail beyond the case of prime powers. Theorem~\ref{thm:main} still reduces to the non-existence of an associated $\Z/p$-equivariant map for a prime~$p$, but generally now $p$ will be much larger than~$q$. The proof method presented here seems to be the first that yields non-trivial upper bounds beyond prime powers, and might turn out to be useful in related contexts. We present the key new idea in Section~\ref{sec:key}.

The continuous generalization of Birch's result is a simple consequence of Theorem~\ref{thm:main}:

\begin{corollary}
\label{cor:main}
	For any drawing of $K_{3q}$ in the plane, where each $3$-cycle is embedded, there are $q$ vertex-disjoint $3$-cycles that surround a common point.
\end{corollary} 

We require that every $3$-cycle needs to be embedded to make sense of the notion of a $3$-cycle surrounding a point: Every $3$-cycle separates $\R^2$ into two regions by the Jordan curve theorem, and the $3$-cycle surrounds every point within the bounded region.

We present some standard results and terminology in Section~\ref{sec:prelim}. Surprisingly, there is a much simpler proof of Theorem~\ref{thm:main} for~${q \le 33}$.  We present it in Section~\ref{sec:simple}. Section~\ref{sec:key} contains the key idea and proofs of Theorem~\ref{thm:main} and Corollary~\ref{cor:main}. The technical verification that an associated configuration space is highly connected is postponed to Section~\ref{sec:construction}.  In Section \ref{sec:final} we present colorful variations of the topological Tverberg theorem beyond prime powers and open problems.

\section{Preliminaries}
\label{sec:prelim}

Here we collect some of the standard language, notation, and results used throughout the manuscript. We refer the reader to Matou\v sek's book~\cite{matousek2003} for an introduction.

We denote the $n$-dimensional simplex $\{x \in \R^{n+1} \ : \ x_i \ge 0 \ \text{and} \ \sum x_i =1\}$ by~$\Delta_n$. A simplicial complex $\Sigma$ is a non-empty collection of sets such that for $\sigma \in \Sigma$ and $\tau \subset \sigma$ we have $\tau \in \Sigma$. All simplicial complexes considered in this manuscript will be finite. By considering every set $\sigma \in \Sigma$ to be a simplex of dimension~$|\sigma|-1$ with the natural identifications, $\Sigma$ is a topological space glued from simplices in a natural way. Our notation does not distinguish between an (abstract) simplicial complex---a collection of finite sets closed under taking subsets---and this geometric realization. For example, as should be obvious from context, a continuous map $f\colon \Sigma \to \R^d$ is defined on the geometric realization of~$\Sigma$. Conversely, we also write $\Delta_n$ for the simplicial complex of all subsets of $[n+1] = \{1,2,\dots, n+1\}$. We refer to $\sigma \in \Sigma$ as a face; its dimension is $|\sigma|-1$. A $0$-dimensional face is called a vertex, a $1$-dimensional face is an edge. A face $\sigma \in \Sigma$ is a maximal face if no proper superset of $\sigma$ is a face of~$\Sigma$.

The join of simplicial complexes $K$ and $L$ is the simplicial complex
$$K * L = \{\sigma\times\{1\} \cup \tau\times\{2\} \ : \ \sigma \in K, \ \tau \in L\}.$$
That is, the faces of $K*L$ are unions of faces of $K$ with faces of~$L$, where we force their vertex sets to be disjoint. The geometric realization of $K*L$ is the join (as topological spaces) of their geometric realizations. The $n$-fold join of $K$ is defined recursively by $K^{*n} = K * K^{*(n-1)}$, where $K^{*1} = K$. A point $x$ in the geometric realization of $K^{*n}$ can thus be represented as an (abstract) convex combination $x = \lambda_1x_1 + \dots + \lambda_n x_n$ for points $x_i$ in~$K$ and $\lambda_i \ge 0$ with $\sum \lambda_i = 1$. 

The homotopical connectivity of a path-connected topological space $X$ will be denoted by~$\conn X$, that is, $\conn X = n$ means that $\pi_j(X, x_0)$ is trivial for $j \le n$ and $\pi_{n+1}(X,x_0)$ is non-trivial for some (and thus any) choice of basepoint~$x_0$. We define $\conn X =-1$ if $X$ is non-empty and not path-connected. We say that $X$ is $n$-connected if $\conn X\ge n$. In particular, a space that is $(n+1)$-connected is also $n$-connected. 

\begin{lemma}[see~{\cite[Prop.~4.4.3]{matousek2003}}]
\label{lem:conn}
	Let $K$ and $L$ be non-empty simplicial complexes. Then $$\conn K*L = \conn K + \conn L + 2.$$
\end{lemma}

We will make repeated use of a simple consequence of the Mayer--Vietoris sequence (and Van Kampen's theorem); see~{\cite[Lemma~10.3]{bjorner1995}}:

\begin{lemma}
\label{lem:mv}
	Let $K$ and $L$ be non-empty simplicial complexes. If $K$ and $L$ are $n$-connected and $K \cap L$ is $(n-1)$-connected then $K \cup L$ is $n$-connected. 
\end{lemma}

A group action of the group $G$ on the space $X$ is called free if $g\cdot x \ne x$ for all non-trivial $g \in G$ and all $x \in X$. If $G$ acts on the spaces $X$ and $Y$ then a continuous map $f\colon X \to Y$ is called $G$-equivariant, or $G$-map, if $f(g\cdot x) = g\cdot f(x)$ for all $g \in G$ and all $x\in X$. If $G$ acts on the simplicial complex $\Sigma$ we will always assume that the action is simplicial, that is, that it maps faces to faces.

\begin{theorem}[Dold, 1983~\cite{dold1983}]
\label{thm:dold}
	Let $G$ be a finite, non-trivial group that acts on the simplicial complex~$\Sigma$ and that acts on~$\R^{n+1}$ by linear maps. Suppose $\Sigma$ is $n$-connected and the action of $G$ restricts to a free action on~$\R^{n+1} \setminus \{0\}$. Then every $G$-map $f\colon \Sigma \to \R^{n+1}$ has a zero, $f(x) = 0$ for some point $x$ in~$\Sigma$.
\end{theorem}

We will need the following often used corollary:

\begin{corollary}
\label{cor:dold}
	Let $p$ be a prime, and let $\Z/p$ act on the simplicial complex~$\Sigma$ and on $(\R^n)^p$ by shifting copies of~$\R^n$. Suppose $\Sigma$ is $[(p-1)n-1]$-connected. Then any $\Z/p$-map $f \colon \Sigma \to (\R^n)^p$ maps some point $x$ in $\Sigma$ to the diagonal $D = \{(y_1, \dots, y_p) \in (\R^n)^p \ : \ y_1 = y_2 = \dots = y_p\}$. 
\end{corollary}

\begin{proof}
	Since $p$ is a prime, the action of $\Z/p$ on the orthogonal complement $D^\perp$ of the diagonal~$D$ is free away from~$0$. Composing $f$ with the orthogonal projection along $D$ onto $D^\perp$ yields an equivariant map $\Sigma \to D^\perp$. This map has a zero since the dimension of $D^\perp$ is~$(p-1)n$.
\end{proof}

\section{Simple proofs for special cases}
\label{sec:simple}

Theorem~\ref{thm:main} is an immediate corollary of the topological Tverberg theorem if $q$ or $q+1$ is a power of a prime. In the latter case we are given a continuous map ${f\colon \Delta_{q(d+1)-1} \to \R^d}$, which we extend continuously to a map ${f^+ \colon \Delta_{q(d+1)} \to \R^d}$ in an arbitrary way by adding a dummy vertex~$v^+$. Since $q+1$ is a prime power, by the topological Tverberg theorem there are $q+1$ pairwise disjoint faces $\sigma_1, \dots, \sigma_{q+1}$ of~$\Delta_{q(d+1)}$ with $${f^+(\sigma_1) \cap \dots \cap f^+(\sigma_{q+1}) \ne \emptyset}.$$ Now simply discard the face that (possibly) contains the dummy vertex~$v^+$, say face~$\sigma_{q+1}$. This leaves $q$ pairwise disjoint faces of $\Delta_{q(d+1)-1}$ with $f(\sigma_1) \cap \dots \cap f(\sigma_q) \ne \emptyset$. For $q+1$ a prime, the optimal colored Tverberg theorem of Blagojevi\'c, Matschke, and Ziegler~\cite{blagojevic2015} gives additional constraints on the faces $\sigma_1, \dots, \sigma_q$.

In this section, we present a particular case of Theorem~\ref{thm:main} that is easier to prove than the general case---the case that $2q+1$ is a prime. This already proves Theorem~\ref{thm:main} for a surprisingly large range of small~$q$. The first positive integer~$q$ that is not a prime power and where $q+1$ is not a prime power and $2q+1$ is not a prime is~$34$. Thus the reduction in this section proves Theorem~\ref{thm:main} in particular for all $q \le 33$.

The main idea of the proof of Theorem~\ref{thm:main} is to consider many copies of the same map~$f$.  This induces a map from a high-dimensional simplex, and with appropriate constraints on the faces we can guarantee that it descends to the desired $q$-fold point of coincidence. For the special case of two copies the required constrained topological Tverberg theorem is known.  The conditions on this Tverberg theorem translate to $2q+1$ being prime.  The result we need is implicit in a paper of Vu\v ci\'c and \v Zivaljevi\'c~\cite{vucic1993}. Denote the vertices of the simplex~$\Delta_n$ by $v_1, v_2, \dots, v_{n+1}$.

\begin{theorem}[Vu\v ci\'c and \v Zivaljevi\'c, 1993~\cite{vucic1993}]
\label{thm:2q+1}
	Let $p$ be an odd prime, and let $d\ge1$ be an integer. For any continuous map $f\colon \Delta_{(p-1)(d+1)} \to \R^d$ there are $p$ pairwise disjoint faces $\sigma_1, \dots, \sigma_p$ such that $f(\sigma_1) \cap \dots \cap f(\sigma_p) \ne \emptyset$ and if $v_{2j-1}$ is in $\sigma_i$ then $v_{2j}$ is in $\sigma_i$ or $\sigma_{i+1}$ for all $j \in \{1,2,\dots, {\frac12(p-1)(d+1)}\}$. Here we consider $\sigma_{p+1}$ to be~$\sigma_1$.
\end{theorem}

Let $f \colon \Delta_{q(d+1)-1} \to \R^d$ be a continuous map, and suppose that $p = 2q+1$ is a prime. Double all vertices of $\Delta_{q(d+1)-1}$ to obtain a map $F$ from the simplex $\Delta_{2q(d+1)-1}$ with twice as many vertices. That is, $F$ linearly interpolates between the map $f$ defined on all vertices $v_2, v_4, \dots, v_{2q(d+1)}$ of $\Delta_{2q(d+1)-1}$ with an even index and the same map~$f$ on the face of $\Delta_{2q(d+1)-1}$ of all odd index vertices $v_1, v_3, \dots, v_{2q(d+1)-1}$. Here we assume that the vertices are ordered in the same way in both copies, that is, if $q\colon \Delta_{2q(d+1)-1} \to \Delta_{q(d+1)-1}$ denotes the linear map defined on vertices by $q(v_{2j-1}) = q(v_{2j}) = v_j$, then $f\circ q = F$.

As before consider the map $F^+\colon \Delta_{2q(d+1)}\to \R^d$ obtained by extending $F$ continuously to another dummy vertex~$v^+ = v_{2q(d+1)+1}$. Now use Theorem~\ref{thm:2q+1} for the map $F^+$ to obtain $2q+1$ pairwise disjoint faces $\sigma_1, \dots, \sigma_{2q+1}$ with $F^+(\sigma_1) \cap \dots \cap F^+(\sigma_{2q+1}) \ne \emptyset$, and with the additional constraint of Theorem~\ref{thm:2q+1}.

By symmetry we may assume that $v^+$ is in~$\sigma_{2q+1}$. We thus get that $${F(\sigma_1) \cap \dots \cap F(\sigma_{2q}) \ne \emptyset}.$$ By definition of~$F$, we get that $f(q(\sigma_1)) \cap \dots \cap f(q(\sigma_{2q})) \ne \emptyset$. Since $q$ identifies pairs of vertices, the faces $q(\sigma_i)$ will generally not be pairwise disjoint. However, by the additional constraint of Theorem~\ref{thm:2q+1} retaining every other face $q(\sigma_1), q(\sigma_3)$, $\dots, q(\sigma_{2q-1})$ yields $q$ pairwise disjoint faces, whose images under $f$ all share a common point.

\section{Key Idea: Sparse, symmetric, highly connected complexes}
\label{sec:key}

Given a continuous map $f\colon\Delta_n \to \R^d$, the condition $f(\sigma_1) \cap \dots \cap f(\sigma_q) \ne \emptyset$ for $q$ pairwise disjoint faces translates to the following equivariant problem: assign a label from the set $[q] = \{1, 2,\dots, q\}$ to each of the $n+1$ vertices of $\Delta_n$ and check if, for the induced partition of the vertices of $\Delta_n$ into $q$ faces, the images of the faces overlap.  This leads us to consider the ``configuration space'' $[q]^{*(n+1)}$, the $(n+1)$-fold join of~$[q]$, where $[q]$ is given the discrete topology. Now the $q$-fold join of~$f$ induces a map $[q]^{*(n+1)} \to (\R^d)^{*q} \subset (\R^{d+1})^q$ that is equivariant with respect to diagonally permuting $[q]$ in the domain and the copies of $\R^{d+1}$ in the codomain. For $q$ a prime power and $n\ge(q-1)(d+1)$ such a map must hit the diagonal $\{(y,\dots,y) \in (\R^{d+1})^q \ : \ y\in \R^{d+1}\}$, which finishes the proof.

In order to prove Theorem~\ref{thm:main}, we generalize this setup. Let $p \ge q$ be a prime, and fix integers $d\ge 1$ and~${n \ge 1}$. Assign to every vertex of $\Delta_n$ a set of labels from~$[p]$.  Each vertex can only receive sets that form faces of a certain simplicial complex $\Sigma$ over~$[p]$, that satisfies the following properties:
\begin{compactenum}[(i)]
	\item\label{it:sym} $\Sigma$ is $\Z/p$-invariant, where the generator $\lambda$ maps $j \mapsto j+1 \mod p$ for all vertices~$j$.
	\item\label{it:sparse} $\Sigma$ has an independent set of size~$q$, that is, $q$ vertices such that no two of them form an edge.
	\item\label{it:conn} The image of every $\Z/p$-equivariant map $\Sigma^{*(n+1)} \to (\R^{d+1})^{p}$ contains a point on the diagonal~$\{(y,\dots,y) \in (\R^{d+1})^p \ : \ y\in \R^{d+1}\}$.
\end{compactenum}

For example, if $p=q$ then property~(\ref{it:sparse}) forces that $\Sigma = [p]$. Property~(\ref{it:conn}) is satisfied whenever $n \ge (p-1)(d+1)$. The larger the prime~$p$ compared to~$q$, the denser the complex $\Sigma$ may be. Our main observation is that for growing~$p$, the denseness of $\Sigma$ required by property~(\ref{it:conn}) may outpace the sparseness imposed by property~(\ref{it:sparse}). In fact, for every prime $p\ge q$ with $p \equiv 1 \mod q$ we will construct a simplicial complex~$\Sigma$ that satisfies properties~(\ref{it:sym}) and~(\ref{it:sparse}), while also being homotopically $(\frac{p}{q} - O(1))$-connected. This implies property~(\ref{it:conn}) for $n\ge q(d+1)$ and sufficiently large~$p$ by Corollary~\ref{cor:dold} and the following:

\begin{lemma}
\label{lem:p-conn}
	Let $c \ge 0$, $d \ge 0$, $q \ge 1$, ${n \ge q(d+1)}$, and $p \ge cq+(c-2)q^2(d+1)$ be integers. Then if $\Sigma$ is a $\left(\frac{p}{q} - c\right)$-connected simplicial complex, its $(n+1)$-fold join $\Sigma^{*(n+1)}$ is ${p(d+1)}$-connected.
\end{lemma}

\begin{proof}
	By Lemma~\ref{lem:conn} $\conn \Sigma^{*(n+1)} \ge (n+1)(\frac{p}{q} - c) + 2n$.
	Now verify that
	\begin{align*}
		&(n+1)\left(\frac{p}{q}-c\right) +2n \ge (q(d+1)+1)\left(\frac{p}{q}-c\right)+2q(d+1) \\
		 = \ &p(d+1)+\frac{p}{q} -c-(c-2)q(d+1) \ge p(d+1),
	\end{align*}
	which completes the proof.
\end{proof}

\begin{remark}
\label{rem:conn}
	The proof of Lemma~\ref{lem:p-conn} shows more generally that if the complex $\Sigma$ is only $\left(\frac{p}{q}-o(p)\right)$-connected, then $\Sigma^{*(n+1)}$ is still ${p(d+1)}$-connected for sufficiently large~$p$.
\end{remark}

\begin{theorem}
\label{thm:Sigma-constraints}
	Let $p$ be a prime, and let $d\ge1$ and $n\ge1$ be integers. Suppose that for these given parameters there is a simplicial complex~$\Sigma$ on~$[p]$ that satisfies properties~(\ref{it:sym}) and~(\ref{it:conn}) above. Then for any continuous map $f\colon \Delta_n \to \R^d$ there are $p$ faces $\sigma_1, \dots, \sigma_p$ of $\Delta_n$ with $f(\sigma_1) \cap \dots \cap f(\sigma_p) \ne \emptyset$ and such that for every vertex $v$ of $\Delta_n$ the set $\{{j \in [p]} \ : \ v \in \sigma_j\}$ is a face of~$\Sigma$.
\end{theorem}

\begin{proof}
	Consider the map $\Phi\colon (\Delta_n)^{*p} \to (\R^{d+1})^p$ defined by
	$$\Phi(\lambda_1x_1 + \dots + \lambda_px_p) = (\lambda_1, \lambda_1f(x_1), \dots, \lambda_p, \lambda_pf(x_p)).$$
	The complex $(\Delta_n)^{*p}$ is isomorphic to~$(\Delta_{p-1})^{*(n+1)}$. Denote the natural isomorphism by $\iota\colon (\Delta_n)^{*p} \to (\Delta_{p-1})^{*(n+1)}$. Thus $\iota^{-1}(\Sigma^{*(n+1)})$ is a $\Z/p$-invariant subcomplex of the domain of~$\Phi$. By property~(\ref{it:conn}) there is a point $x = \lambda_1x_1 + \dots + \lambda_px_p \in (\Delta_n)^{*p}$ with $\iota(x) \in \Sigma^{*(n+1)}$ and $\Phi(x) = (y, \dots, y)$, or equivalently, 
	$$\lambda_1 = \lambda_2 = \dots = \lambda_p \quad \text{and} \quad \lambda_1f(x_1) = \lambda_2f(x_2) = \dots = \lambda_pf(x_p).$$
	This implies $\lambda_j = \frac1p$ for all~$j$, and thus $f(x_1) = f(x_2) = \dots = f(x_p)$.
	
	\begin{figure}[h]
\centerline{\includegraphics[scale=1]{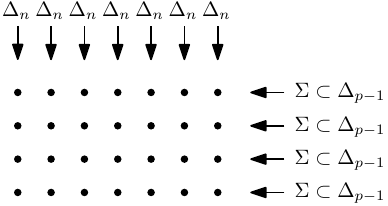}}
\caption{The simplex on vertex set~$[n+1] \times [p]$ is both the join of columns $(\Delta_n)^{*p}$ as well as the join of rows~$(\Delta_{p-1})^{*(n+1)}$. Here $\Sigma$ is a subcomplex of each row. The isomorphism $\iota$ reverses the roles of rows and columns.}
\end{figure}
	
	Now let $\sigma_j$ be the inclusion-minimal face of $\Delta_n$ that contains~$x_j$. Then $f(\sigma_1) \cap \dots \cap f(\sigma_p) \ne \emptyset$. Moreover $\iota(\sigma_1 * \dots * \sigma_p)$ is a face of~$\Sigma^{*(n+1)}$, but this means that each vertex $v$ may only appear in faces $\sigma_j$ whose indices form a face of~$\Sigma$.
\end{proof}

\begin{corollary}
\label{cor:Sigma-constraints}
	Let $p$ be a prime, and let $d\ge1$, $q\ge2$, and $n\ge1$ be integers. Suppose that for these given parameters there is a simplicial complex~$\Sigma$ on~$[p]$ that satisfies properties~(\ref{it:sym}), (\ref{it:sparse}), and~(\ref{it:conn}) above.  Then for any continuous map $f\colon \Delta_{n-1} \to \R^d$ there are $q$ pairwise disjoint faces $\sigma_1, \dots, \sigma_q$ of $\Delta_{n-1}$ with $f(\sigma_1) \cap \dots \cap f(\sigma_q) \ne \emptyset$.
\end{corollary}

\begin{proof}
	Let us extend $\Delta_{n-1}$ to $\Delta_{n}$ by adding a dummy vertex $v^+$.  We can extend $f$ as well, turning it into a map $f:\Delta_n \to \rr^d$. By Theorem~\ref{thm:Sigma-constraints} there are $p$ faces $\sigma_1, \dots, \sigma_p$ of $\Delta_n$ with $f(\sigma_1) \cap \dots \cap f(\sigma_p) \ne \emptyset$ and such that for every vertex $v$ of $\Delta_n$ the set $\{j \in [p] \ : \ v \in \sigma_j\}$ is a face of~$\Sigma$.  Let $\sigma^+ = \{j \in [p] \ : \ v^+ \in \sigma_j\}$. Now let $I \subset [p]$ be an independent set of size~$q$.  We claim there exists an $m \in [p]$ such that $\lambda^m I \cap \sigma^+ = \emptyset$, where $\lambda$ is a generator of $\mathbb{Z}/p$.  If $\sigma^+ = \emptyset$ the claim is true, so we may assume $|\sigma^+| \ge 1$.   The set $\lambda^m I$ is independent for each $m$.  If our claim fails to be true, consider the set
	\[
	P = \left\{(i,m) : i \in I, m \in [p], \lambda^m i \in \sigma^+ \right\}.
	\] 
	
	For each $m$, we know $|\lambda^m I \cap \sigma^+|\ge 1$ since it is not empty.  We also know $|\lambda^m I \cap \sigma^+|\le 1$ since it is the intersection of a complete set with an independent set in $\Sigma$.  Therefore, $|P| = p$.  However, for each $i \in I$, we know there are exactly $|\sigma^+|$ values of $m$ such that $(i,m) \in P$.  Since $p>q$ is prime, it cannot be equal to $|\sigma^+|q$.  Therefore, there exists some $m$ such that $\lambda^m I \cap \sigma^+ = \emptyset$.
	
	 The faces $\sigma_j$ with $j \in \lambda^m I$ are pairwise disjoint, and none contains the dummy vertex~$v^+$. After possibly renumbering the $\sigma_j$ we thus have that $f(\sigma_1) \cap \dots \cap f(\sigma_q) \ne \emptyset$ for pairwise disjoint faces $\sigma_1, \dots, \sigma_q$ of $\Delta_{n-1}$.
\end{proof}

To prove Theorem~\ref{thm:main} for given $q\ge2$ and $d\ge1$, we have to find a prime~$p$ and a $\Z/p$-invariant complex $\Sigma$ on~$[p]$ with an independent set of size~$q$ that is $(\frac{p}{q}-O(1))$-connected. We will construct such a complex $\Sigma$ in the next section, thus proving our main result:

\begin{proof}[Proof of Theorem~\ref{thm:main}]
	Let the integers $q\ge2$ and $d\ge1$ be given. There are infinitely many primes $p$ of the form $p = (a+1)q+1$ by Dirichlet's theorem~\cite{dirichlet1837}. By Theorem~\ref{thm:cyclic-conn} for each such prime $p$ there is a simplicial complex $C^a_p$ that satisfies properties~(\ref{it:sym}) and~(\ref{it:sparse}) and is $(\frac{p}{q}-4)$-connected. For $n = q(d+1)$ the $(n+1)$-fold join $(C_p^a)^{*(n+1)}$ is $p(d+1)$-connected for sufficiently large~$p$ by Lemma~\ref{lem:p-conn}. Thus $C_p^a$ satisfies property~(\ref{it:conn}) as well by Corollary~\ref{cor:dold}. Thus for any continuous $f \colon \Delta_{n-1} \to \R^d$ there are $q$ pairwise disjoint faces whose images all share a common point by Corollary~\ref{cor:Sigma-constraints}.
\end{proof}

\begin{proof}[Proof of Corollary~\ref{cor:main}]
	Let $f\colon K_{3q} \to \R^2$ be a continuous map such that every $3$-cycle is embedded. The complete graph $K_{3q}$ is the $1$-skeleton of~$\Delta_{3q-1}$. Let $\sigma$ be a triangle of~$\Delta_{3q-1}$. By the Jordan curve theorem $f(\partial\sigma)$ bounds a disk $D$ in~$\R^2$. Continuously extend $f$ to $\sigma$ such that $f(\sigma) \subset D$. After extending $f$ onto every triangle of~$\Delta_{3q-1}$, continuously extend it to all of~$\Delta_{3q-1}$. By Theorem~\ref{thm:main} there are $q$ pairwise disjoint faces $\sigma_1, \dots, \sigma_q$ with $f(\sigma_1) \cap \dots \cap f(\sigma_q) \ne \emptyset$. By the methods presented by Sch\"oneborn and Ziegler~\cite[Sec.~2]{schoneborn2005} the $\sigma_i$ can be chosen such that $\dim\sigma_i \le 2$, which completes the proof.
\end{proof}

We remark that our construction of sparse, $\Z/p$-symmetric, highly connected complexes in Section~\ref{sec:construction} is optimal, as otherwise by the same reasoning as above one could prove Tverberg-type results that are too strong to be true:

\begin{theorem}
	Let $q\ge2$ be an integer. Then for any sufficiently large prime $p$ the maximal size of an independent set in any $\Z/p$-invariant complex on~$[p]$ that is $(\frac{p}{q}-o(p))$-connected is~$q$.
\end{theorem}

\begin{proof}
	Suppose $\Sigma$ is a simplicial complex on~$[p]$ that is $\Z/p$-invariant, $(\frac{p}{q}-o(p))$-connected, and has an independent set of size~${q+1}$. By Lemma~\ref{lem:p-conn} and Remark~\ref{rem:conn} for $p$ sufficiently large $\Sigma$ satisfies property~(\ref{it:conn}) above for $n=q(d+1)$ for any given~$d\ge1$. By the reasoning in the proof of Corollary~\ref{cor:Sigma-constraints} for any continuous map $f\colon  \Delta_{n-1} \to \R^d$ there are $q+1$ pairwise disjoint faces $\sigma_1, \dots, \sigma_{q+1}$ of $\Delta_{n-1}$ with $f(\sigma_1) \cap \dots \cap f(\sigma_{q+1}) \ne \emptyset$. However, the bound of Tverberg's theorem is known to be optimal. That is, an affine map $f$ that maps vertices into sufficiently generic position yields a contradiction.
\end{proof}

\section{Construction of simplicial complexes}
\label{sec:construction}

In this section we construct the required $\Z/p$-invariant complex $\Sigma$ that has an independent set of size~$q$ and is almost $(p/q)$-connected. Throughout this section the integer~$q\ge2$ will be fixed. Suppose $p = (a+1)q+1$ is a prime for some integer $a\ge 1$. Observe that there are infinitely many such integers~$a$ by Dirichlet's theorem~\cite{dirichlet1837}. A set $\sigma \subset \Z/p$ is \emph{$q$-stable} if for any two $i,j \in \sigma$ the difference $i-j$ is not in~$\{1,2,\dots,q-1\}$, that is, the cyclic gap between any two elements of $\sigma$ is at least~$q$. Let $C_p$ be the simplicial complex of all $q$-stable subsets of~$\Z/p$, and let $C_p^a$ be the subcomplex of $q$-stable sets that can be extended to a $q$-stable set of size at least~$a$, that is,
$$C_p^a = \{\sigma \subset \Z/p \ : \ \sigma \ \text{is} \ q\text{-stable and there is a} \ \tau \in C_p \ \text{with} \ \sigma \subset \tau \ \text{and} \ |\tau| \ge a\}.$$

Clearly $C_p^a$ is $\Z/p$-invariant, and any $q$ consecutive vertices form an independent set. The rest of this section is devoted to checking that $C_p^a$ is highly connected:

\begin{theorem}
\label{thm:cyclic-conn}
	For $p = (a+1)q+1$ the complex $C_p^a$ is $(a-2)$-connected. 
\end{theorem}

We remark that $a-2 = \frac{p-1}{q}-3 \ge \frac{p}{q} -4$, so the theorem above finishes the proof of the main result. To prove Theorem~\ref{thm:cyclic-conn} we will first break the cyclic symmetry and study subcomplexes of $C_p^a$ obtained by restricting to $p-q+1$ successive vertices.

The extendability condition seems necessary to guarantee high connectedness.  For instance, if $q=2$, the complex $C_r$ is the independence complex of a cycle of length $r$.  Its connectedness is $\left\lfloor (r-1)/3\right\rfloor - 1$, as its homotopy type was determined by Kozlov~\cite{Kozlov:1999ik} to be a wedge of spheres of dimension  $\left\lfloor (r-1)/3\right\rfloor$.  For general~$q$, the complex $C_r$ has maximal faces which are roughly $(r/(2q-1))$-dimensional, which may cause the connectedness to drop.

The complexes $C_p^a$ have also appeared in the study of minimal manifold triangulations: K\"uhnel and Lassmann~\cite{kuhnel1996} study $C_p^{a+1}$ for $p = (a+1)q+1$ and certain generalizations since they are triangulations of disk bundles over the circle. In particular, $C_p^{a+1}$ is not highly connected, while $C_p^a$ is.

\subsection{Linear complexes and their connectedness}

Fix an integer~$r \ge q$. We define the simplicial complex $L_r$ as the family of sets $\sigma \subset [r]$ such that if $j \in \sigma$, then none of $j+1, \ldots, j+q-1$ are in $\sigma$. (Here addition is \textbf{not} modulo~$r$.) In other words, any two elements of $\sigma$ have a difference greater than or equal to~$q$.  In analogy to their cyclic counterparts, define
$$L^a_r = \{\sigma \in L_r : \mbox{ there exists } \tau \in L_n \mbox{ such that } |\tau| \ge a, \sigma \subset \tau \}.$$

In other words, the maximal faces of $L^a_r$ have at least $a$ vertices. This makes the complexes $L_r^a$ and $C_p^a$ highly connected, whereas $L_r$ and $C_p$ are not as highly connected. 

If $r \ge aq$, then the vertex set of $L^a_n$ is~$[r]$.  However, if $(a-1)q+1 \le n < aq$, not all elements of $[r]$ appear as vertices.  For example, $L^a_{(a-1)q+1}$ has a single face $\sigma = \{1,q+1, \ldots, (a-1)q+1\}$ and no other vertices.

For a simplicial complex $T$ with vertex set in~$[r]$ and an integer~$s$, we define $T+s$ as the complex obtained by translating each face of $T$ by $s$ units, with vertex set in $\{1+s,2+s,\dots,r+s\}$.

\begin{claim}
\label{claim:L}
	The complex $L^a_r$ is either empty or $(a-2)$-connected.
\end{claim}

\begin{proof}
	We proceed by induction on~$a$.  If $a=1$, the result is clear.  Now suppose $a \ge 2$ and $r \ge (a-1)q+1$, so $L^a_r$ is not empty.  Every maximal face of $L^a_r$ has exactly one of the vertices $\{1,2,\dots, q\}$.  Let $A_i$ be the subcomplex of faces in $L^a_r$ that contain or can be extended to contain vertex~$i$.  Notice that $A_i$ is contractible or empty (if there are no faces containing vertex~$i$). Then,
	\[
	L^a_r = A_1 \cup \dots \cup A_q.
	\]
	
	For $k=1,\dots, q$, let $B_k = A_1 \cup \dots \cup A_k$.  We show by induction on $k$ that $B_k$ is $(a-2)$-connected.  Since $L^a_r$ is not empty, then $A_1$ is not empty, so $B_1 = A_1$, which is contractible.  If we know that $B_k$ is $(a-2)$-connected for some $k<q$, then notice $B_{k+1} = B_k \cup A_{k+1}$.  The first complex is $(a-2)$-connected and the second is contractible.  We only have to show that $B_{k} \cap A_{k+1}$ is $(a-3)$-connected.  This intersection is the set of faces that can be extended to a face of at least $a-1$ vertices among $\{q+k+1, \dots, r\}$, so
	\[
	B_k \cap A_{k+1} = L^{a-1}_{r-k-q} + (k+q).
	\]

By induction, we know this complex is either empty or $(a-3)$-connected.  If this intersection is $(a-3)$-connected, then $B_{k+1}$ is $(a-2)$-connected.  However, if this intersection is empty that means that $A_{k+1}$ must be empty, so $B_{k+1} = B_k$.  In either case $B_{k+1}$ is $(a-2)$-connected.
\end{proof}

A simple corollary of this claim is a lower bound on the connectedness of $L_r$ for any $r$.

\begin{corollary}
	The complex $L_r$ is $\left(\left\lfloor \frac{r}{2q-1}\right\rfloor-2\right)$-connected.
\end{corollary}

\begin{proof}
	The corollary follows immediately from the fact that, for $a = \left\lfloor \frac{r}{2q-1}\right\rfloor$, $L_r = L^a_r$.  The containment $L^a_r \subset L_r$ is true by definition of $L^a_r$.  A maximal face of $L_r$ cannot leave a gap of length $2q-1$ between two consecutive elements, or it could be extended.  Therefore maximal faces of $L_r$ have at least $a$ elements, implying $L_r \subset L^a_r$.
\end{proof}

Next, we need to characterize all faces of the complex $L^{a-1}_r$ for small values~$r$. If $r \le (a-2)q$, then $L^{a-1}_r$ is empty. If $(a-2)q +1 \le r \le (a-1)q$, then the faces of $L^{a-1}_r$ can be fully characterized.  A maximal face has exactly $a-1$ vertices in $\{1,\ldots, r\}$.  These vertices can be partitioned into $a-1$ blocks $I_1, \ldots, I_{a-1}$, where for $j=1, \ldots, a-2$
\[
I_j = \{(j-1)q+1, \ldots, jq\}
\]
and $I_{a-1} = \{(a-2)q+1, \ldots, r\}$.

\begin{figure}[h]
\centerline{\includegraphics[scale=1]{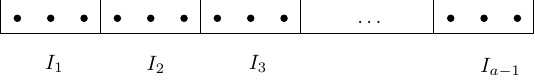}}
\caption{Blocks in $L^{a-1}_r$}
\end{figure}

A maximal face of $L^{a-1}_r$ must have exactly one element of each block.  If the element of the $j$-th block is $(j-1)q+k_j$, we know that
\[
1 \le k_1 \le \dots \le k_{a-1} \le r-(a-2)q \le q.
\]

We need a particular complex $T^{a-1}_m \subset L^{a-1}_m$ for $m= (a-2)q+k+2$, where $1 \le k \le q-1$.  We define $T^{a-1}_m$ by its maximal faces. A maximal face $\tau$ of $L^{a-1}_m$ is a maximal face of $T^{a-1}_m$ if~$\tau$
\begin{compactitem}
	\item does not contain vertex~$1$, or
	\item does not contain any of the last $k$ vertices.
\end{compactitem}

That is, $T^{a-1}_m$ is the subcomplex of $L^{a-1}_m$ such that its maximal faces have sequences $(k_1, \ldots, k_{a-1})$ that satisfy either $k_1 \ge 2$ or $2 \ge k_{a-1}$.

\begin{claim}
\label{claim:T}
 For $m= (a-2)q+k+2$ the complex $T^{a-1}_m$ is $(a-3)$-connected.
\end{claim}

\begin{proof}
	The complex $T^{a-1}_m$ is the union of two sets of maximal faces in~$L^{a-1}_m$: Those maximal faces whose sequences satisfy $k_1 \ge 2$ and those whose sequences satisfy ${k_{a-1} \le 2}$. The former set of faces form the complex~$L^{a-1}_{m-1}+1$, while the latter produce~$L^{a-1}_{m-k}$. Since $m-k = (a-2)q+2$ both of these complexes are non-empty, and thus $(a-3)$-connected by Claim~\ref{claim:L}.  Moreover, since $m-1 \le (a-1)q$ for the range of $k$ we are considering, the characterization of their faces above still holds.  Their intersection consists of all those faces of~$L^{a-1}_m$ contained in a face of size $a-1$ with $k_1 \ge 2$ and $k_{a-1} \le 2$. Since the sequence of $k_i$ is increasing, we have $k_1 = k_2 =\dots = k_{a-1} = 2$, which implies that the intersection is a simplex, and thus contractible. Thus $T^{a-1}_m$ is $(a-3)$-connected.
\end{proof}

\begin{claim}\label{claim-uniones}
	For $1 \le k \le q-1$ and $r = aq+2$ the complex
	\[
	L^a_{r-1} \cup \left(L^{a-1}_{r+k-2q} + (q-1)\right)
	\]
	is $(a-3)$-connected.
\end{claim}

\begin{proof}
	Since each of $L^a_{r-1}$ and $\left(L^{a-1}_{r+k-2q} + (q-1)\right)$ is $(a-3)$-connected, we only have to check that their intersection $K = L^a_{r-1} \cap \left(L^{a-1}_{r+k-2q} + (q-1)\right)$ is $(a-4)$-connected.
	
	The vertices of the complex $K$ are contained in $\{q,q+1,\dots,r+k-q-1\}$. We first treat the case that $k \le q-2$. Then $r+k-2q \le (a-1)q$ and so no face of $L^{a-1}_{r+k-2q} + (q-1)$ can have more than $a-1$ elements. Thus the maximal faces of $K$ are exactly those maximal faces of $L^{a-1}_{r+k-2q} + (q-1)$ that can be extended by one of the vertices $\{1,\dots, q-1\}$ or $\{r+k-q,\dots, r-1\}$ for it to be in~$L^a_{r-1}$. We claim that $K$ is~${T^{a-1}_{r+k-2q} + (q-1)}$, which is $(a-3)$-connected by Claim~\ref{claim:T}.
	
	We first show that $K \subset {T^{a-1}_{r+k-2q} + (q-1)}$. Let $\sigma$ be a maximal face of~$K$. If it can be extended to a face of size $a$ in $L^a_{r-1}$ by a vertex in $\{1,\dots, q-1\}$, then $\sigma$ cannot contain the vertex~$q$. Thus $\sigma$ is in~${T^{a-1}_{r+k-2q} + (q-1)}$. If on the other hand $\sigma$ may be extended by a vertex in~$\{r+k-q,\dots, r-1\}$, then it cannot contain vertices~$r-q, r-q+1, \dots, r-q+k-1$. These are precisely the last $k$ vertices of the vertex set of~${L^{a-1}_{r+k-2q} + (q-1)}$ and thus $\sigma$ is in~${T^{a-1}_{r+k-2q} + (q-1)}$.
	
	Now if $\sigma$ is a maximal face of~${T^{a-1}_{r+k-2q} + (q-1)}$, then it has size~${a-1}$. If $\sigma$ does not contain vertex~$q$, then it is contained in the face $\sigma \cup \{1\}$ of $L^a_{r-1}$ and thus in~$K$. If $\sigma$ does not contain any of the last $k$ vertices of the vertex set of~${T^{a-1}_{r+k-2q} + (q-1)}$ then it is contained in the face $\sigma \cup \{r-1\}$ of $L^a_{r-1}$ and thus in~$K$. This shows that $K$ is indeed~${T^{a-1}_{r+k-2q} + (q-1)}$ and thus $(a-3)$-connected.
	
	Lastly, we deal with the case $k = q-1$. The analysis is the same as above with the only difference that now a maximal face of $K$ may have $a$ vertices. Since $r+k-2q = (a-1)q+1$, there is a unique such face $\sigma^* = \{q,2q,\ldots, aq\}$, and $K$ is the complex 
	$$\sigma^* \cup \left(T^{a-1}_{r+k-2q} + (q-1)\right).$$
	Here we treat $\sigma^*$ as a simplicial complex and tacitly include all of its subsets as faces. This exhibits $K$ as the union of a simplex, which is contractible, and an $(a-3)$-connected complex. The intersection of $\sigma^*$ and $T^{a-1}_{r+k-2q} + (q-1)$ consists of all those faces of~$\sigma^*$ that do not use both $q$ and~$aq$. Thus their intersection is a bipyramid (suspension of an $(a-2)$-simplex), which is contractible. This completes the proof.
\end{proof}

\subsection{Circular arcs complexes and their connectedness.}

We now study the complexes~$C_p^a$. There is a striking difference when studying the connectedness of $L^a_r$ and that of~$C^a_p$.  The complex $L^a_r$ is either empty or highly connected.  With $C^a_p$ this fails to hold.  For example, for $p=aq$, the complex $C^a_p$ consists of $q$ pairwise disjoint simplexes, so it is non-empty and disconnected.  We need more vertices to be able to guarantee high connectedness.

First, let us relate the linear complexes and the circular complexes:
\begin{claim}
\label{claim:L-to-C}
For every $a$ and $r$ we have
\[
C^a_{r+q-1} = \bigcup_{j=0}^{q-1} \left(L^a_r + j\right)
\]
	
\end{claim}

\begin{proof}
	Let $\sigma$ be a non-empty maximal face of $C^a_{r+q-1}$.  Let $s \in \{1,2,\dots, r+q-1\}$ be the largest element of~$\sigma$.  If $s \le r$, then $\sigma$ is in~$L^a_r$.  Otherwise, let $j = s-r$.  Note that $1 \le j \le q-1$.  We claim that $\sigma$ is in $L^a_r +j$.  For this we only need to verify that $\sigma \subset \{j+1, \dots, r+j\}$.  By the choice of~$s$, $\sigma$ does not contain any larger elements.  By $q$-stability, $\sigma$ cannot contain any of $s+1, \dots, s+q-1$, all modulo $r+q-1$.  Notice that $s+q-1$ modulo $r+q-1$ is~$j$, and so the smallest element vertex of $\sigma$ is larger than~$j$.
	
	Now let $\sigma$ be a maximal face of~$L^a_r+j$.  The cyclic gap between the smallest possible element~$j+1$ and the largest possible element~$r+j$ is $q$, and thus $\sigma$ is a face of~$C^a_{r+q-1}$.
\end{proof}

The following claim shows Theorem~\ref{thm:cyclic-conn}.

\begin{claim}\label{claim-circular-connected}
	For $r=aq+2$, the complex $C^{a}_{r+q-1}$ is $(a-2)$-connected.
\end{claim}

\begin{proof}
	By Claim~\ref{claim:L-to-C}, we know that
	\[
C^a_{r+q-1} = \bigcup_{j=0}^{q-1} \left(L^a_r + j\right)
\]
	
	For $k=0,\dots, q-1$, consider the complex 
	\[
M_k = \bigcup_{j=0}^{k} \left(L^a_r+ j\right)
\]	
	We prove by induction on $k$ that $M_k$ is $(a-2)$-connected.
	
	\noindent
	\underline{Base of induction}. $M_0 = L^a_r$, which is $(a-2)$-connected.
	
	\noindent
	\underline{Inductive step}.  Assume that $M_{k-1}$ is $(a-2)$-connected for some $1 \le k\le q-1$ and we want to show that $M_k$ is also $(a-2)$-connected.  We can write $M_k$ as
	\[
	M_k = M_{k-1} \cup (L^a_r+k)
	\]
	The first complex is $(a-2)$-connected by induction, and the second is $(a-2)$-connected by Claim~\ref{claim:L}, so we only need to look at their intersection.  The only vertices we can take in the intersection are contained in $\{k+1,k+2,\dots, r+k-1\}$.  We claim that $M_{k-1}\cap \left(L^a_r + k\right) = \left( L^a_{r-1} + k \right) \cup \left(L^{a-1}_{r+k-2q}+q\right)$.
	
	Let $\sigma$ be a face of~$M_{k-1}\cap \left(L^a_r + k\right)$. Thus $\sigma$ is contained in a maximal face $\tau_1$ of~$\left(L^a_r + k\right)$ and in a maximal face $\tau_2$ of~$\left(L^a_r + j\right)$ for some $j \in \{0,1,\dots, k-1\}$. If $r+k\notin \tau_1$ then $\tau_1$ and thus $\sigma$ is a face of~$\left(L^a_{r-1} + k\right)$. If $r+k \in \tau_1$, then $\tau_1 \cap \tau_2$ has at most $a-1$ vertices and $r+k-q+1, r+k-q+2, \dots, r+k-1$ are not in~$\tau_1\cap \tau_2$. If $\tau_2$ does not contain one of the first $k$ vertices, we could choose $\tau_1 = \tau_2$, which implies that $\sigma$ is a face of~$L^a_{r-1}+k$; a case that we have already treated above. But if $\tau_2$ does contain one of the first $k$ vertices, then $\tau_1\cap \tau_2$ does not contain any vertex among $\{1,2,\dots,q\}$. Thus $\tau_1 \cap \tau_2$ is a face of~$\left(L^{a-1}_{r+k-2q}+q\right)$, and so $\sigma$ is too.
	
	If $\sigma$ is a face of~$(L^a_{r-1}+k)$, then it is a face of~$(L^a_r+k)$ and of~$(L^a_r+(k-1))$ and thus in $M_{k-1}\cap \left(L^a_r + k\right)$. If $\sigma$ is a maximal face of~$\left(L^{a-1}_{r+k-2q}+q\right)$, then $\sigma \cup \{r+k\}$ is a face of~$L_r^a+k$, and $\sigma \cup \{1\}$ is a face of~$L^a_r$. Thus $\sigma$ is a face of~$M_{k-1}\cap \left(L^a_r + k\right)$. 
	
	We now show that $\left( L^a_{r-1} + k \right) \cup \left(L^{a-1}_{r+k-2q}+q\right)$ is $(a-3)$-connected.  Notice that among the vertices $\{1,2,\dots, r+k-1\}$ the map $x \mapsto r+k-x$ induces an isomorphism of the complexes
	\begin{align*}
		\left(L^a_{r-1}+k \right) & \to L^a_{r-1} \\
		\left(L^{a-1}_{r+k-2q} + q\right) & \to \left( L^{a-1}_{r+k-2q} + (q-1)\right)
	\end{align*} 
	
	Therefore, $M_{k-1} \cap \left(L^a_r + k\right) \cong L^a_{r-1} \cup \left( L^{a-1}_{r+k-2q} + (q-1)\right)$, which by Claim~\ref{claim-uniones} is $(a-3)$-connected.
\end{proof}

\section{Final remarks}\label{sec:final}

We hope that this manuscript is only the starting point for many variants of Theorem~\ref{thm:main} beyond prime powers. Combining our main result with the work of Blagojevi\'c, the first author, and Ziegler~\cite{blagojevic2014} yields some such variants ``for free,'' for example the following colorful version follows immediately from combining Theorem~\ref{thm:main} with~\cite{blagojevic2014}. It extends the colorful Tverberg theorem of \v Zivaljevi\'c and Vre\'cica~\cite{zivaljevic1992} beyond prime powers at the expense of allowing for a higher-dimensional simplex.

\begin{theorem}
	Let $q\ge2$, $c\ge 0$, and $d\ge1$ be integers. Let $C_1, \dots, C_c$ be sets of vertices of~$\Delta_{q(d+c+1)-1}$ each of size at most~${2q-1}$. Let $f\colon \Delta_{q(d+c+1)-1} \to \R^d$ be a continuous map. Then there are $q$ pairwise disjoint faces $\sigma_1, \dots, \sigma_q$ of $\Delta_{q(d+c+1)-1}$ such that $f(\sigma_1) \cap \dots \cap f(\sigma_q) \ne \emptyset$ and each $\sigma_i$ has at most one vertex in each~$C_j$.
\end{theorem}

If we wish to reduce the number of points in each color class $C_j$, the reader may verify that it is sufficient to have $|C_j| \ge q + \frac{q(d+1)}{c}$, and the domain to be a simplex with $c$ times as many vertices.  Therefore, as $c \to q(d+1)$, the number of elements required in each color class approaches $q+1$.  This gives us a continuous Tverberg type-theorem similar in parameters to the ``equal coefficients'' colorful Tverberg \cite{blagojevic2014, Soberon:2015cl} that works beyond prime powers. 

We provide a precise statement below. The proof is the same as~\cite[Thm.~8.1]{blagojevic2014}, now using the new cases beyond prime powers of Theorem~\ref{thm:main}. Let $C$ be a set of vertices of the simplex~$\Delta_n$. Let $x$ and $y$ be two points in~$\Delta_n$, and thus they can be written as convex combination of the vertices $v_1,\dots,v_{n+1}$, say $x = \sum \lambda_iv_i$ and $y = \sum \mu_iv_i$. We say that $x$ and $y$ have equal barycentric coordinates with respect to~$C$ if $\sum_{v_i \in C} \lambda_i = \sum_{v_i \in C} \mu_i$.

\begin{theorem}
	Let $q\ge2$ and $d\ge1$ be integers. Let $C_1 \sqcup \dots \sqcup C_{q(d+1)}$ be a partition of the vertex set of~$\Delta_{q(q+1)(d+1)-1}$ each of size~${q+1}$. Let $f\colon \Delta_{q(q+1)(d+1)-1} \to \R^d$ be a continuous map. Then there are $x_1, \dots, x_q$ in $q$ pairwise disjoint faces of $\Delta_{q(q+1)(d+1)-1}$ such that $f(x_1) = \dots = f(x_q)$ and all $x_i$ have the same barycentric coordinates with respect to each~$C_j$. 
\end{theorem}

It remains open if a colorful version of Theorem~\ref{thm:main} is true.  This is a known conjecture by B\'ar\'any and Larman~\cite{Barany:1992tx}, that has not been settled even for affine maps.  It has been verified for $q+1$ a prime by Blagojevi\'c, Matschke, and Ziegler~\cite{Blagojevic:2011vh, blagojevic2015}.

\begin{conjecture}
	Let $q \ge 2$, $d \ge 1$ be integers.  Let $C_1, \ldots, C_{d+1}$ be pairwise disjoint sets of vertices of of~$\Delta_{q(d+1)-1}$ each of size $q$. Let $f\colon \Delta_{q(d+1)-1} \to \R^d$ be a continuous map. Then there are $q$ pairwise disjoint faces $\sigma_1, \dots, \sigma_q$ of $\Delta_{q(d+1)-1}$ such that $f(\sigma_1) \cap \dots \cap f(\sigma_q) \ne \emptyset$ and each $\sigma_i$ has at most one vertex in each~$C_j$.
\end{conjecture}

\section*{Acknowledgements}

\noindent
Thanks to G\"unter Ziegler for helpful comments on a draft of this manuscript.

\bibliographystyle{amsalpha}

\end{document}